\documentclass[preprint,ECP]{ejpecp}

\usepackage[T1]{fontenc}
\usepackage[utf8]{inputenc}

\usepackage{stmaryrd,enumitem}

\SHORTTITLE{Overdamped limit  for non-equilibrium Langevin diffusions}

\TITLE{Overdamped limit at stationarity for non-equilibrium Langevin diffusions}

\AUTHORS{%
 Pierre Monmarch{\'e}\footnote{LJLL and LCT, Sorbonne Université, Paris, France.
    \EMAIL{pierre.monmarche@sorbonne-universite.fr}}
  \and 
  Mouad Ramil\footnote{CERMICS, Ecole des Ponts, Marne-la-Vallée, France. \EMAIL{ramil.mouad@gmail.com} }}


\AMSSUBJ{60J60 ; 82C31} 

\SUBMITTED{October 4, 2021} 
\ACCEPTED{?} 

\VOLUME{0}
\YEAR{2020}
\PAPERNUM{0}
\DOI{10.1214/YY-TN}

\ABSTRACT{
In this note, we establish that the stationary distribution of a possibly non-equilibrium Langevin diffusion converges, as the damping parameter goes to infinity (or equivalently in the Smoluchowski-Kramers vanishing mass limit), toward a tensor product of the stationary distribution of the corresponding overdamped process and of a Gaussian distribution.}

 \newcommand{\po}{\left(}
\newcommand{\pf}{\right)}

\newcommand{\cco}{\llbracket}
\newcommand{\ccf}{\rrbracket}
\newcommand{\R}{\mathbb R} 
\newcommand{\T}{\mathbb T} 
\newcommand{\Z}{\mathbb Z} 
 
\newcommand{\dd}{\mathrm{d}}
\newcommand{\na}{\nabla}

\begin{document}



\section{Introduction}

Consider the Langevin diffusion  on $\mathbb T^d\times\R^d$ ($\T=\R/\Z$) solution of 
\begin{equation}\label{eq:Langevin}
\left\{\begin{array}{rcl}
\dd X_t & = & Y_t\dd t \\
\dd Y_t & = & F(X_t) \dd t - \gamma Y_t\dd t + \sqrt{2\gamma} \Sigma \dd B_t,
\end{array}\right.
\end{equation}
with a standard $d$-dimensional Brownian motion $(B_t)_{t\geqslant0}$, a force field $F\in \mathcal C^1(\T^d,\R^d)$, a damping parameter $\gamma>0$ and a constant symmetric positive definite diffusion matrix $\Sigma$. For all $\gamma>0$, this Markov process admits a unique invariant measure $\mu_\gamma$. Besides, as  $\gamma\rightarrow +\infty$,  $(X_{\gamma t})_{t\geqslant 0}$ converges to the overdamped Langevin diffusion on $\T^d$, solution of
\begin{equation}\label{eq:overdamped}
\dd Z_t \ = \ F(Z_t)\dd t + \sqrt 2\Sigma\dd B_t\,,
\end{equation}
which admits a unique invariant measure $\mu_O$. These facts are well-known,  see e.g. \cite{LelievreStoltz}. Up to a rescaling,  this overdamped limit is equivalent to the so-called Smoluchowski-Kramers limit, where a mass parameter is sent to zero. 
In the so-called equilibrium case, in which $F=-\na U$ is a conservative force ($U\in\mathcal C^2(\T^d)$) and $\Sigma=I$, $\mu_\gamma$ and $\mu_O$ are explicit, and more precisely, for all $\gamma>0$, $\mu_\gamma=\mu_O\otimes g$ with
\[ \mu_O(\dd x) = \frac{e^{-U(x)}}{\int_{\mathbb T^d} e^{-U(w)}\dd w}\dd x\,,\qquad g(\dd y) = \frac{e^{-|y|^2/2}}{(2\pi)^{d/2}}\dd y\,.\]
Here and in the rest of this work, $\T^d$ is endowed with the Lebesgue measure $\dd x$ normalized so that its total mass is $1$.  
In the general non-equilibrium case, however, there is no reason for any of this to be true, namely $\mu_O$ and $\mu_\gamma$ are not explicit, $\mu_\gamma$ does depend on $\gamma$, it is not a tensor product and its position marginal is not $\mu_O$ (see Section~\ref{sec:explicit}). The main result of this note is  that, nevertheless, these properties are recovered in the overdamped limit, i.e. that $\mu_\gamma$ converges as $\gamma\rightarrow +\infty$ toward $\mu_O\otimes g_\Sigma$ where $g_\Sigma$ is the centered Gaussian distribution on $\R^d$ with covariance matrix $\Sigma^2$. The convergence is quantitative, and stated in term of  Wasserstein distance. The Wasserstein distance between two probability laws $\mu,\nu$ on a metric space $(E,\mathrm{dist})$ is defined as
 \[\mathcal W(\mu,\nu) = \inf_{\pi\in \Gamma(\mu,\nu)} \int_{E^2} \mathrm{dist}(y,z)\pi(\dd z,\dd y)\,,\]
 where $\Gamma(\mu,\nu)$ is the set of probability measures on $E^2$ with marginals $\mu$ and $\nu$. Our result is the following.

\begin{theorem}\label{thm:main}
There exists $C>0$, which depends on $F,\Sigma$ and $d$, such that, for all $\gamma\geqslant 2$,
\[\mathcal W(\mu_\gamma,\mu_O \otimes g_\Sigma) \leqslant C \frac{\sqrt{\log\gamma}}{\gamma}\,.\]
\end{theorem} 

\begin{remark}
We consider a position in the compact torus since this is often the case considered in previous works on non-equilibrium Langevin processes, see e.g. \cite{IacobucciOllaStoltz}. Moreover, standard molecular dynamics simulations use  periodic boundary conditions \cite{LMS,LelievreStoltz,RousselStoltz}. However, it is straightforward to adapt the proof to a position on $\R^d$ under the additional assumptions that there exist $R,\lambda>0$ such that $F(x)\cdot x \leqslant - \lambda |x|^2$ for all $|x|\geqslant R$ and that the Jacobian matrix of $F$ is bounded.  
\end{remark}

This note is organized as follows. We finish this introduction by discussing applications and previous related works. In Section~\ref{sec:explicit} we characterize the cases where, as in the equilibrium case, position and velocity are independent under $\mu_\gamma$, or $\mu_\gamma$ is independent from $\gamma$. Section~\ref{sec:proof} is devoted to the proof of Theorem~\ref{thm:main}.

\bigskip
Non-equilibrium Langevin processes appear in various situations where a physical system undergoes an external force or a gradient of temperature, see e.g. \cite{StoltzProc,Fourier,LetiziaOlla,RousselStoltz,EVANS19901} or \cite[Chapter 5]{LelievreStoltz} and references within. In these cases, macroscopic quantities of interest such as transport coefficients (e.g.  mobility,  thermal conductivity or shear viscosity) can be expressed in term of expectations with respect to the stationary state of the non-equilibrium process, see e.g.  \cite{EVANS19901} for details.

More precisely, with $\Sigma=I$, the question of the linear response of a system at equilibrium to a non-conservative perturbation leads to the study of forces of the form $F(x)=-\na U(x)+\tau \tilde F(x)$ where $\tilde F$ is non-conservative and $\tau$ is small. For instance, \cite[Theorem 5.2]{LelievreStoltz} gives an asymptotic expansion of $\mu_\gamma$ as $\tau$ vanishes (at a fixed $\gamma$). Similarly, the long-time convergence of non-equilibrium Langevin processes is studied in \cite{IacobucciOllaStoltz} in this case, for $\tau$ small enough. Our result is not restricted to this perturbative regime.

Chains of oscillators such as those studied in \cite{Fourier,LetiziaOlla} and references within, which are models for thermal conductivity, correspond to a case where $F(x)=-\na U(x)$ for some potential $U$ but $\Sigma$ is not an homothety. Hence, the Jacobian matrix of $\Sigma^2 F$ is not symmetric, which shows that $\Sigma^2 F$ cannot be a gradient (which means this is a non-equilibrium process, cf. Section~\ref{sec:explicit}).

For more details and a general discussion on non-equilibrium Langevin processes from a mathematical perspective, we refer to \cite[Chapter 5]{LelievreStoltz}.

\bigskip

To the best of our knowledge, Theorem~\ref{thm:main}, or even the weak convergence of the first marginal of $\mu_\gamma$ toward $\mu_O$, is new. A similar result as been proven by one of the authors in \cite{Ram} but, rather than the stationary measure of the process, it concerns the quasi-stationary measure of the process killed when it exists a given domain. This is a different and more difficult problem, besides the convergence in \cite{Ram} is not quantitative. In infinite dimension, a result similar to Theorem~\ref{thm:main} (although not quantitative and concerning only the convergence of the position marginal), with some similarity in the proof, is established in \cite{Cerrai} for stochastic damped wave equations.

Distinct but somewhat related results are the overdamped limit of the diffusion coefficient proven in \cite{HairerPavliotis,LMS} in dimension 1 and of the large deviation quasi-potential in \cite{ChenFreidlin}. Indeed, these two works study the overdamped limit of a stationary quantity in some asymptotic regime, either $\tau\rightarrow 0$ (with $F=-\na U + \tau \tilde F$) or $\Sigma\rightarrow 0$.

\section{Particular cases}\label{sec:explicit}

In this section, we investigate under which conditions on $F$ and $\Sigma$ each of the following properties holds true: a) $\mu_\gamma$ is independent from $\gamma$, b) $\mu_\gamma$ is the tensor product of two probability laws respectively on $\T^d$ and $\R^d$, and c) the first marginal of $\mu_\gamma$ is $\mu_O$. Before addressing the general case, let us highlight the  following two simple cases:
\begin{itemize}
    \item \emph{The equilibrium case.} If $F(x)=-\Sigma^2 \na U(x)$ for some $U\in\mathcal C^2(\T^d)$   then the properties a), b) and c) hold since it is easily checked that
    \[\mu_\gamma(\dd x,\dd y) = \frac{e^{-U(x)-|\Sigma^{-1}y|^2/2}}{(2\pi)^{d/2}\mathrm{det}(\Sigma)\int_{\T^d} e^{-U(w)}\dd w} \dd x \dd y\]
    for all $\gamma>0$. Considering on $\Sigma^{-1}\T^d$ the potential $\tilde U(z)=U(\Sigma z)$, the process $(Z,W)=(\Sigma^{-1}X,\Sigma^{-1}Y)$ on $(\Sigma^{-1}\T^d,\R^d)$ solves
    \[\left\{\begin{array}{rcl}
         \dd Z_t &=& W_t\dd t   \\
         \dd W_t &= & -\na \tilde U(Z_t)\dd t - \gamma W_t \dd t + \sqrt 2 \dd B_t\,. 
    \end{array}\right.\]
    \item \emph{The space-homogeneous case.} If $F(x)=\eta$ for all $x\in\T^d$ for some constant $\eta\in\R^d$, then 
    \[\mu_\gamma(\dd x,\dd y) = \frac{e^{-|\Sigma^{-1}(y-\eta/\gamma)|^2/2}}{(2\pi)^{d/2}\mathrm{det}(\Sigma)} \dd x \dd y\]
    for all $\gamma>0$. In particular, the properties b) and c) hold but a) holds iff $\eta=0$ (which corresponds to the equilibrium case with a constant $U$). Moreover, since the Wasserstein distance between two Gaussian distributions with the same variance is given by the distance between their means (see e.g.  \cite[Proposition 7]{W2Gauss}),
    \[\mathcal W\po \mu_\gamma,\mu_O\otimes g_\Sigma\pf = \frac{|\Sigma^{-1}\eta|}{\gamma}.\] 
\end{itemize}

Let us state a characterization of the cases where a) or b) hold.

\begin{proposition}\label{prop:particular}
First, the two following properties are equivalent.
\begin{enumerate}[label=(\roman*)]
    \item There exist $\gamma_1>\gamma_2>0$ such that $\mu_{\gamma_1}=\mu_{\gamma_2}$.
    \item There exists $U\in\mathcal C^2(\T^d)$ such that $F=-\Sigma^2 \na U$.
\end{enumerate}
Second, the two following properties are equivalent.
\begin{enumerate}[label=(\roman*)]
  \setcounter{enumi}{2}
    \item There exists $\gamma>0$ such that $\mu_\gamma$ is the tensor product of two probability laws respectively on $\T^d$ and $\R^d$.
    \item There exists $\eta\in\R^d$ and $U\in\mathcal C^2(\T^d)$ such that, for all $x\in\T^d$, $ F(x)=-\Sigma^2 \na U(x)+\eta$ and $\na U(x)\cdot \eta=0$.
\end{enumerate}
Moreover, if \textit{(iv)} holds, then
\[\mu_O(\dd x) = \frac{e^{-U(x)}}{\int_{\T^d} e^{-U(w)}\dd w} \dd x\qquad \text{and for all } \gamma>0,\qquad \mu_\gamma = \mu_O \otimes \frac{e^{-|\Sigma^{-1}(y-\eta/\gamma)|^2/2}}{(2\pi)^{d/2}\mathrm{det}(\Sigma)} \dd y. \]
\end{proposition}

Before proving this result, let us discuss the condition $\na U\cdot \eta=0$ in $(iv)$, which is equivalent to say that $U(x+t\eta)=U(x)$ for all $x\in\T^d,t\in\R$. Since $U$ is $\mathbb Z$-periodic, depending on $\eta$, this may be more or less restrictive. For instance:
\begin{itemize}
    \item If $\R\eta$ is dense in $\T^d$ (which is in particular the case if $d=1$ and $\eta\neq 0$), then necessarily $U$ is constant, hence so is $F$.
    \item If $\eta = \tau \mathbf{e}$ with $\tau \in\R$ and $\mathbf{e}$ a vector of the canonical basis (say, the first one), then the condition is satisfied as soon as $U$ is periodic function of $x^{\neq 1}=(x_j)_{j\in\cco 2,d\ccf}$. If moreover, for instance, $\Sigma=I$, the process can be decomposed as two independent parts, a position homogenous coordinate and the other coordinates at equilibrium, namely
    \begin{equation}\label{eq:decouple}
    \left\{\begin{array}{rcl}
         \dd X^1_t &=&  Y^1_t \dd t \\
         \dd Y^1_t &=& \tau  \dd t -\gamma Y^1_t\dd t + \sqrt{2\gamma} \dd B_t^1  \\
         \dd X^{\neq 1}_t &=&  Y^{\neq 1}_t \dd t \\
         \dd Y^{\neq 1}_t &=& -\na_{\neq 1} U(X^{\neq 1}) \dd t -\gamma Y^{\neq 1}_t\dd t + \sqrt{2\gamma} \dd B_t^{\neq 1}  \,. 
    \end{array}\right.
    \end{equation}
   \end{itemize}
In any cases, assuming $(iv)$ with $\eta\neq 0$, we see that the process $(Z,W)=(Q\Sigma^{-1}X,Q\Sigma^{-1} W)$ on $Q\Sigma^{-1}\T^d\times\R^d$ where $Q$ is an orthonomal matrix that sends $\eta/|\eta|$ to $\mathbf{e}$, solves an equation of the form \eqref{eq:decouple}. In other words, apart from the equilibrium case, the only cases where $\mu_\gamma$ is of a tensor form for some $\gamma>0$ is a combination of a one-dimensional space homogeneous system with an independent $d-1$ dimensional equilibrium system.

\begin{proof}

 The Langevin process being hypoelliptic and irreducible, its invariant measure admits a smooth positive density with respect to the Lebesgue measure, which we write $e^{-H(x,y)}$ ($H$ possibly depends on $\gamma$ but we do not explicit this dependency to simplify notation). The invariance of $\mu_\gamma$ is equivalent to
\[\forall f\in\mathcal C^\infty(\T^d\times\R^d)\,,\qquad \int_{\T^d \times\R^d} L f(x,y) e^{-H(x,y)}\dd x \dd y =0\,, \]
where, using the notation $\Sigma^2 :\na^2_y = \sum_{i,j=1}^d (\Sigma^2)_{i,j} \partial_{y_i} \partial_{y_j}$,
\begin{equation}\label{eq:generateur}
    L = y\cdot \nabla_x   + (F(x)-\gamma y)\cdot \na_y + \gamma \Sigma^2 :\na^2_y
\end{equation}
is the generator of the process \eqref{eq:Langevin}. Integrating by parts, this is equivalent to
\begin{multline}\label{eq:A}
   \forall (x,y)\in\T^d\times\R^d\,,\qquad   
0 = y \cdot \na_x H(x,y) + \po F(x)-\gamma y\pf\cdot  \na_y H(x,y) \\
+ \gamma d  + \gamma  | \Sigma \na_y H(x,y)|^2   - \gamma \Sigma^2 :\na^2_y  H(x,y).
\end{multline}
The implications $(ii)\rightarrow (i)$ and $(iv)\rightarrow(iii)$ immediately follow from the fact that \eqref{eq:A} holds with $H(x,y)=U(x)+|\Sigma^{-1}(y-\eta/\gamma)|^2/2$ if $F=-\Sigma^2\na U +\eta$ and $\na U\cdot\eta=0$.

Let us prove that $(i)\rightarrow (ii)$. Assuming $(i)$ and applying \eqref{eq:A} for $\gamma\in\{\gamma_1,\gamma_2\}$ yields
\[\forall (x,y)\in\T^d\times\R^d\qquad \left\{\begin{array}{rcl}
y \cdot \na_x H(x,y) +  F(x)\cdot   \na_y H(x,y) & = & 0 \\ 
- y\cdot  \na_y H(x,y) + d  + | \Sigma \na_y H(x,y)|^2  -  \Sigma^2 :\na^2_y  H(x,y) & =& 0\,.
\end{array}
\right.
\]
The second equation is equivalent to say that for all $x\in\T^d$, $y\mapsto e^{-H(x,y)}/\int_{\R^d} e^{-H(x,z)}\dd z$ is the invariant probability density of the Ornstein-Uhlenbeck process $\dd Y = -Y\dd t + \sqrt 2\Sigma\dd B_t$, in other words is the $d$-dimensional Gaussian distribution with mean $0$ and covariance matrix $\Sigma^{2}$. This means that $H(x,y)=U(x)+| \Sigma^{-1} y|^2/2$ for some $U\in\mathcal C^1(\T^d)$. Plugging this in the first equation yields $ F(x)=-\Sigma^2\na U(x)$ for all $x\in\T^d$.

We now turn to the proof of the implication $(iii)\rightarrow (iv)$. Let $\gamma$ be such that $\mu_\gamma$ is a tensor product, namely $H(x,y)=U(x)+V(y)$ for some $U\in\mathcal C^1(\T^d)$, $V\in\mathcal C^1(\R^d)$. Let $x_*\in\T^d$ be a critical point of $U$. Applying \eqref{eq:A} for $x=x_*$ reads
\[\forall y\in\R^d\qquad  \po F(x_*)-\gamma y\pf\cdot  \na V(y) + \gamma d  + \gamma |\Sigma \na V(y)|^2  - \gamma \Sigma:\na^2 V(y) = 0\,,\]
in other words $e^{-V}$ is an invariant measure for the Ornstein-Uhlenbeck process
\[\dd Y_t = (F(x_*)-\gamma Y_t)\dd t + \sqrt {2\gamma} \Sigma \dd B_t\,.\]
Hence, $V(y)= |\Sigma^{-1}(y-F(x_*)/\gamma)|^2/2$. 
Plugging this back in \eqref{eq:A} yields
\[\forall (x,y)\in\T^d\times\R^d\,,\qquad y \cdot \na_x U(x)+ \po F(x)-F(x_*)\pf\cdot \Sigma^{-2}  \po y-F(x_*)/\gamma\pf    = 0 \,, \]
which is equivalent to
\[\forall x\in\T^d\,,\qquad F(x) -  F(x_*) = -\Sigma^2 \na U(x) \,,\qquad (F(x)-F(x_*))\cdot\Sigma^{-2} F(x_*)=0\,,\]
and thus $\na U(x)\cdot F(x_*)=0$ for all $x\in\T^d$.

Finally, we prove the last statement of the proposition as follows. The expression of $\mu_\gamma$ under $(iv)$ is proven by checking  \eqref{eq:A}. Since the position marginal doesn't depend on $\gamma$, the fact that it is necessarily equal to $\mu_O$  follows from letting $\gamma\rightarrow +\infty$ in Theorem~\ref{thm:main} (it can also be checked directly on the overdamped stationary equation).

\end{proof}

\section{Proof of Theorem~\ref{thm:main}}\label{sec:proof}

Let $(P_t^\gamma)_{t\geqslant0}$ be the Markov semi-group associated with the Langevin process  \eqref{eq:Langevin} with friction $\gamma$ and $(\tilde P_t)_{t\geqslant 0}$ be the semi-group associated with the overdamped process \eqref{eq:overdamped}, namely
\[P_t^\gamma f(x,y) = \mathbb E_{(x,y)} \po f(X_t,Y_t)\pf\,,\qquad \tilde P_t f(z) = \mathbb E_z \po f(Z_t)\pf\]
for all bounded measurable functions $f$ respectively on $\T^d\times\R^d$ and $\T^d$. If $\nu$ is a probability on $\T^d$ then $\nu \tilde P_t$ is the law of $Z_t$ if $Z_0\sim \nu$, and similarly for $P_t^\gamma$.

\begin{lemma}\label{lem:Eberle}
There exists $C_1,\rho>0$, which depends only on $F,d$ and $\Sigma$, such that, for all probability measures $\mu,\nu$ on $\T^d$ and all $t\geqslant 0$,
 \[\mathcal W(\mu \tilde P_t,\nu \tilde P_t) \leqslant C_1 e^{-\rho t} \mathcal W(\mu,\nu)\,.\]
\end{lemma}
\begin{proof}
This is proven in \cite{Eberle}, although the latter is written in $\R^d$ and not $\T^d$. The adaptation is straightforward.
\end{proof}

The following is a simple fact on Wasserstein distances.

\begin{lemma}\label{lem:marginalWasserstein}
Let $(E_i,r_i)$ for $i=1,2$ be two metric spaces and endow $E_1\times E_2$ with a distance $r$ such that $r((x,y),(z,w))\geqslant r_1(x,z)$ for all $(x,y),(z,w)\in E_1\times E_2$. Let $\mu,\nu$ be two probability measures on $E_1\times E_2$ and $\mu^1,\nu^1$ be their marginal laws on $E_1$. Then
\[\mathcal W(\mu^1,\nu^1) \leqslant \mathcal W(\mu,\nu).\]
\end{lemma}
\begin{proof}
Let $((X,Y),(Z,W))\sim \pi \in \Gamma(\mu,\nu)$ be a coupling of $\mu$ and $\nu$. Then the law of $(X,Z)$ is in $\Gamma(\mu^1,\nu^1)$, and thus
\[\mathcal W(\mu^1,\nu^1) \leqslant \mathbb E(r_1(X,Z)) \leqslant \mathbb E(r((X,Y),(Z,W))).\]
Taking the infimum over $\Gamma(\mu,\nu)$ concludes the proof.
\end{proof}

Next, we state a moment bound, uniform in $\gamma\geqslant 1$.

\begin{lemma}\label{lem:moment}
For all $\gamma>0$,
\[ \mathbb E_{\mu_\gamma}\po |Y|^2\pf \leqslant 2 \mathrm{Tr}(\Sigma^2) + \frac{\|F\|_\infty^2}{\gamma^2} \,.\]
\end{lemma}
\begin{proof}
Since $\mu_\gamma$ is a stationary measure for the Langevin process, considering $f(x,v)=|v|^2/2$ and the generator $L$ given by \eqref{eq:generateur}, 
\begin{align*}
0 = \int_{\T^d\times\R^d} Lf(x,v)\mu_\gamma(\dd x,\dd v) & \leqslant \int_{\T^d\times\R^d} \po \|F\|_\infty |y| - \gamma |y|^2 +\gamma \mathrm{Tr}(\Sigma^2)\pf \mu_\gamma(\dd x,\dd v)     \\
& \leqslant - \frac\gamma 2 \int_{\T^d\times\R^d}  |y|^2  \mu_\gamma(\dd x,\dd v)     + \gamma \mathrm{Tr}(\Sigma^2) + \frac{\|F\|_\infty^2}{2\gamma}\,,
\end{align*}
which concludes the proof.
\end{proof}

The next result is based on some arguments from \cite{Ram}.

\begin{proposition}\label{prop:couple}
For all $t>0$, there exists $C_2>0$, which depends only on $t,F,d$ and $\Sigma$, such that, for all $\gamma\geqslant 2$ and all probability measure $\nu$ on $\T^d\times\R^d$, denoting by $\nu^1$ the marginal of $\nu$ on $\T^d$,
\[\mathcal{W}\left(\nu P^{\gamma}_{\gamma t},\nu^1\Tilde{P}_{t}\otimes g_\Sigma\right)\leq C_2\frac{\sqrt{\log(\gamma)}}{\gamma}\left(1+\mathbb{E}_{\nu}\po\vert Y\vert\pf\right).\]
\end{proposition}
Before proving this proposition, let us first  conclude the proof of our main result.

\begin{proof}[Proof of Theorem~\ref{thm:main}]
Notice that in the statement of the theorem, the distance used on $\T^d\times\R^d$ is implicit, namely it is any distance equivalent to the standard $r((x,y),(z,w)) = \mathrm{dist}(x,z)+|y-w|$, where $\mathrm{dist}$ is the image of the Euclidean norm $|\cdot|$ on the torus. We prove the theorem with this distance.

Considering $C_1 $ and $\rho$ as in Lemma~\ref{lem:Eberle}, let $t=\log(2C_1)/\rho$. Using that $\mu_\gamma$ and $\mu_O$ are invariant respectively for $P_{\gamma t}^\gamma$ and $\tilde P_t$, and denoting by $\mu_\gamma^1$ the marginal law of $\mu_\gamma$ on $\mathbb T^d$, we get
\begin{eqnarray*}
\mathcal W \po \mu_\gamma,\mu_O\otimes g_\Sigma\pf & = & \mathcal W \po \mu_\gamma P_{\gamma t}^\gamma,\po  \mu_O \tilde P_t\pf\otimes g_\Sigma\pf\\
& \leqslant & \mathcal W \po \mu_\gamma P_{\gamma t}^\gamma, \po\mu_\gamma^1 \tilde P_{t}\pf\otimes g_\Sigma\pf + \mathcal W \po \po \mu_\gamma^1 \tilde P_{t}\pf \otimes g_\Sigma,\po \mu_O \tilde P_t\pf\otimes g_\Sigma\pf\\
& \leqslant & \mathcal W \po \mu_\gamma P_{\gamma t}^\gamma, \po \mu_\gamma^1 \tilde P_{t}\pf\otimes g_\Sigma\pf + C_1 e^{-\rho  t}\mathcal W \po \mu_\gamma^1 ,  \mu_O \pf\\
& \leqslant & \mathcal W \po \mu_\gamma P_{\gamma t}^\gamma, \po\mu_\gamma^1 \tilde P_{t}\pf\otimes g_\Sigma\pf + \frac12\mathcal W \po \mu_\gamma,\mu_O\otimes g_\Sigma\pf\,,
\end{eqnarray*}
where we used Lemma~\ref{lem:marginalWasserstein} for the last inequality. Proposition~\ref{prop:couple} then yields
\[\mathcal W \po \mu_\gamma,\mu_O\otimes g_\Sigma\pf \leqslant 2 C_2\frac{\sqrt{\log(\gamma)}}{\gamma}\left(1+\mathbb{E}_{\mu_\gamma}\po\vert Y\vert\pf\right)\] 
and Lemma~\ref{lem:moment} concludes the proof.
\end{proof}

\begin{proof}[Proof of Proposition~\ref{prop:couple}]
We are going to construct an explicit coupling introduced in~\cite{Ram}. The time $t>0$ is fixed. Let $(X_0,Y_0)$ be a random variable with law $\nu$ and $(B_s)_{s\geqslant 0}$ be a standard Brownian motion independent from $(X_0,Y_0)$. For $\gamma>0$, let $(X_s^{(\gamma)},Y_s^{(\gamma)})_{s\geqslant 0}$ be the solution  of \eqref{eq:Langevin} driven by $(B_s)_{s\geqslant 0}$ with initial condition $(X_0,Y_0)$,  define the Brownian motion $(B^{(\gamma)}_s)_{s\geq0}=(B_{\gamma s}/\sqrt{\gamma})_{s\geq0}$ and let $(\overline{X}^{(\gamma)}_s)_{s\geq0}$ be the overdamped process starting from $X_0$ and solving \eqref{eq:overdamped} where the Brownian motion is replaced by the $\gamma$-dependent Brownian motion $(B^{(\gamma)}_s)_{s\geq0}$. Let also
\begin{equation}\label{def A_t gamma}
    \forall s\geq0,\qquad A_{s}^{(\gamma)}:= \sqrt 2 \Sigma\mathrm{e}^{-\gamma^2 s} \int_0^{\gamma s}\mathrm{e}^{\gamma r} \mathrm{d}B_r.
\end{equation}
Then, for all $s\geq0$,
\begin{equation}\label{ecriture Y_gamma t}
    Y^{(\gamma)}_{\gamma s}= \mathrm{e}^{-\gamma^2 s}Y_0 +\gamma \mathrm{e}^{-\gamma^2 s}\int_0^s\mathrm{e}^{\gamma^2 r} F(X^{(\gamma)}_{\gamma r}) \mathrm{d}r+A_s^{(\gamma)}.
\end{equation}
Moreover, let $h^{(\gamma)}_t:[0,t]\mapsto\mathbb{R}$ and the process $(Z_{s,t}^{(\gamma)})_{s\in[0,t]}$ be defined as follow:
\begin{equation}\label{expr h}
    \forall s\in[0,t],\qquad h^{(\gamma)}_t(s):=\frac{2}{\gamma}\frac{\mathrm{e}^{- \gamma^2 (t-s)}-\mathrm{e}^{- \gamma^2 t}}{1-\mathrm{e}^{-2 \gamma^2 t}},\qquad Z_{s,t}^{(\gamma)}:=\Sigma B^{(\gamma)}_s-h^{(\gamma)}_t(s) A_t^{(\gamma)}.
\end{equation}  

Let $(\mathcal{F}^{(\gamma),Z}_s)_{s \in [0,t]}$ be the natural filtration of $(Z_{s,t}^{(\gamma)})_{s\in[0,t]}$. Using a Itô’s fixed point argument, see for instance~\cite[Thm 2.9 p. 289]{Karatzas}, the following equation       
\begin{equation*}
  \mathrm{d}W^{(\gamma)}_s=F(W^{(\gamma)}_s) \mathrm{d}s+\sqrt 2 \mathrm{d}Z_{s,t}^{(\gamma)}\,,\qquad W_0^{(\gamma)}=X_0
\end{equation*}
possesses a unique strong solution on $[0,t]$ which is adapted to $(\mathcal{F}^{(\gamma),Z}_s)_{s \in [0,t]}$. Moreover, it is shown in~\cite[Lemma 3.1]{Ram} that the process $(W^{(\gamma)}_s)_{s\in[0,t]}$ is independent of the random variable $A_t^{(\gamma)}$. Denoting by $ \mathcal L(Q)$ the law of a random variable $Q$, we bound
\begin{multline*}
    \mathcal{W}\left(\nu P^{\gamma}_{\gamma t},\nu^1\Tilde{P}_{t}\otimes g_\Sigma\right)  \leq
    \mathcal{W}\left(\nu P^{\gamma}_{\gamma t}, \mathcal L(W_t^{(\gamma)},A_t^{(\gamma)})\right)
    + \mathcal{W}\left(\mathcal L(W_t^{(\gamma)})\otimes \mathcal L(A_t^{(\gamma)}),\nu^1\Tilde{P}_{t}\otimes g_\Sigma\right)\\
     \leq \mathbb E \po \mathrm{dist}\po X_{\gamma t}^{(\gamma)},W_t^{(\gamma)}\pf + |Y_{\gamma t}^{(\gamma)}-A_t^{(\gamma)}|\pf + \mathbb E \po\mathrm{dist}\po \overline{X}_{t}^{(\gamma)},W_t^{(\gamma)} \pf \pf + \mathcal W\po \mathcal L(A_t^{(\gamma)}), g_\Sigma\pf .
\end{multline*}
It remains to bound each of these terms. First, it follows from~\cite[(i) Lemma 3.2]{Ram} that there exists a constant $C>0$ such that for all $\gamma\geqslant1$, 
$$\mathbb{E}\left[\mathrm{dist}\po  X^{(\gamma)}_{\gamma t},W^{(\gamma)}_t\pf \right]\leq\frac{C}{\gamma}\left(1+\mathbb{E}\po|Y_0|\pf+ \sqrt{\log(1+\gamma^2t)}\right)\mathrm{e}^{Ct}.$$
Second, from~\eqref{ecriture Y_gamma t},
$$\mathbb{E}\left[\left\vert Y^{(\gamma)}_{\gamma t}-A_t^{(\gamma)}\right\vert\right]\leq \mathbb{E}\po |Y_0|\pf \mathrm{e}^{-\gamma^2 t}+\frac{\Vert F\Vert_\infty}{\gamma}.$$
Third, it follows from~\cite[(ii) Lemma 3.2]{Ram} that there exists a constant $C>0$ such that for all $\gamma\geqslant 1$,
$$\mathbb{E}\po \mathrm{dist}\po  W^{(\gamma)}_t,\overline{X}^{(\gamma)}_t\pf\pf\leq\frac{C}{\gamma}\mathrm{e}^{Ct}.$$
Finally, since $\mathcal L(A_t^{(\gamma)})$ is a centered Gaussian distribution with variance  $\Sigma^2(1-\mathrm{e}^{-2 \gamma^2 t})$, the formula for the $L^2$ Wasserstein distance (which is larger than the $L^1$ distance $\mathcal W$) between Gaussian distributions (see \cite{W2Gauss}) yields
$$\mathcal{W}\left(\mathcal{L}(A_t^{(\gamma)}),\mu_\Sigma\right)\leq\sqrt{\mathrm{Tr}\left(\Sigma^2\right)}\left(1-\sqrt{1-\mathrm{e}^{-2\gamma^2t}}\right).$$
Summing the four last inequalities concludes the proof.
\end{proof}

\textbf{Acknowledgements:} P. Monmarché acknowledges financial support by the  French ANR grants EFI (ANR-17-CE40-0030) and SWIDIMS (ANR-20-CE40-0022) and by the European Research Council (ERC) under the European Union's Horizon 2020 research and innovation program (grant agreement No 810367),
project EMC2.  Mouad Ramil is currently a postdoc at CEA-DAM in Arpajon (France). The authors thank Gabriel Stoltz for fruitfull discussions.




\bibliographystyle{amsplain}
\bibliography{biblio}






\end{document}